\documentclass{amsart}

\usepackage{amsmath}
\usepackage{amsfonts}
\usepackage{amssymb,enumerate}
\usepackage{amsthm}
\usepackage{cite}
\usepackage{comment}
\usepackage[all]{xy}
\usepackage{hyperref}

\theoremstyle{definition}
\newtheorem{lem}{Lemma}[section]
\newtheorem{cor}[lem]{Corollary}

\newtheorem{thm}[lem]{Theorem}

\newtheorem{Defn}[lem]{Definition}
\newtheorem{Ex}[lem]{Example}
\newtheorem{Question}[lem]{Question}
\newtheorem{Property}[lem]{Property}
\newtheorem{Properties}[lem]{Properties}
\newtheorem{Discussion}[lem]{Remark}
\newtheorem{Construction}[lem]{Construction}
\newtheorem{Notation}[lem]{Notation}
\newtheorem{Fact}[lem]{Fact}
\newtheorem{Notationdefinition}[lem]{Definition/Notation}
\newtheorem{Remarkdefinition}[lem]{Remark/Definition}

\newtheorem{Subprops}{}[lem]
\newtheorem{Para}[lem]{}
\newtheorem{Exer}[lem]{Exercise}
\newtheorem{Exerc}{Exercise}

\newenvironment{defn}{\begin{Defn}\rm}{\end{Defn}}
\newenvironment{ex}{\begin{Ex}\rm}{\end{Ex}}

\renewcommand{\geq}{\geqslant}
\renewcommand{\leq}{\leqslant}

\numberwithin{equation}{section}

\begin{document}
\author[J. G. Boynton]{Jason Greene Boynton}
\address{Department of Mathematics\\
	North Dakota State University\\
	Fargo, ND 58108}
\email[J. G. Coykendall]{jason.boynton@ndsu.edu}
\author[J. Coykendall]{Jim Coykendall}
\address{Department of Mathematics\\
	North Dakota State University\\
	Fargo, ND 58108}
\email[J. Coykendall]{jim.coykendall@ndsu.edu}

\title[Graph of divisibility]{On the graph of divisibility of an integral domain}
\begin{abstract}
It is well-known that the factorization properties of a domain are reflected
in the structure of its group of divisibility. The main theme of this paper
is to introduce a topological/graph-theoretic point of view to the current
understanding of factorization in integral domains. \ We also show that
connectedness properties in the graph and topological space give rise to a
generalization of atomicity.
\end{abstract}

\maketitle

\section{Introduction}

Let $D$ be an integral domain with field of fractions $K.$ \ Then, the group
of divisibility $G(D)$ is defined to be the partially ordered additive group
of principal fractional ideals with $aD\leq bD$ if and only if $aD\supseteq
bD.$ \ If $K^{\times }$ is the multiplicative group of $K$ and if $U(D)$ is
the group of units of $D$, then $G(D)$ is order isomorphic to the quotient
group $K^{\times }/U(D)$ with the ordering $aU(D)\leq bU(D)$ if and only if $%
\frac{b}{a}\in D.$ \ 

It is well-known that the factorization properties of a domain are reflected
in the structure of its group of divisibility. \ For example, an integral
domain is a unique factorization domain if and only if its group of
divisibility is a direct sum of copies of $%
%TCIMACRO{\U{2124} }%
%BeginExpansion
\mathbb{Z}
%EndExpansion
$ equipped with the usual product order. \ It is also true that the group of
divisibility reflects more than just factorization properties of a domain. \
Indeed, it is not hard to check that a domain is a valuation domain if and
only if its group of divisibility is totally ordered. \ We refer the
interested reader to \cite{Mott} for an excellent survey of material
regarding the group of divisibility.

In 1968, Cohn introduced the notion of an atomic integral domain in \cite%
{Cohn}. \ These are the domains in which every nonzero nonunit admits a
finite factorization into irreducible elements. \ For several years, it was
believed to be the case that atomicity in an integral domain was equivalent
to the ascending chain condition on principal ideals (ACCP). \ However, in
1974, Anne Grams demonstrated that an atomic domain need not satisfy ACCP in 
\cite{Gr}. \ Grams was able to understand the subtle difference between
atomicity and ACCP using the group of divisibility. \ Ten years later, Zaks
added two more examples of an atomic domain without ACCP in \cite{Z}. \
However, examples of atomic domains without ACCP are still relatively scarce.

The main theme of this work is to introduce a topological/graph-theoretic
point of view to the current understanding of factorization in integral
domains. \ That is, we find a graphical representation of the group of
divisibility in order to detect various well-studied factorization
properties of an integral domain. \ The contents of this paper is organized
as follows. \ In Section 2, we recall a topological structure that is
naturally associated to a partially ordered set. \ In addition, we make the
relevant graph-theoretic definitions needed in the sequel. \ In Section 3,
we introduce the \textit{graph of divisibility} of an integral domain and
show that this graph detects the standard factorization properties studied
in \cite{AAZ}. \ In Section 4, we examine the connectedness properties of
the graph of divisibility using some elementary topology to do so. \ In
Section 5, we will see that a connected graph of divisibility gives rise to
a generalized atomicity. \ We also provide some examples in order to
illustrate these notions.

\section{Some Definitions and Background}

In this section, we make some relevant definitions from graph theory and
topology that will be used throughout. \ We refer the reader to \cite{Arenas}
for a survey of known results about the Alexandrov topology. \ 

\begin{defn}
Let $(X,\tau )$ be a topological space with neighborhood base $\mathcal{U}%
(x)=\{U\in \tau :x\in U\}.$

\begin{enumerate}
\item $(X,\tau )$ is called an \textit{Alexandrov space} if arbitrary
intersections of open sets remain open. \ 

\item For every $x$ in an Alexandrov space $X$, we set $M(x)=\cap _{U\in 
\mathcal{U}(x)}U.$ \ The set $M(x)$ is called the \textit{minimal open set}
containing $x.$
\end{enumerate}
\end{defn}

\begin{thm}
Let $(X,\tau )$ be an Alexandrov space.

\begin{enumerate}
\item The collection of minimal open sets $\mathcal{N}=\{M(x):x\in X\}$ is a
basis for the space $(X,\tau ).$ \ 

\item $(X,\tau )$ is a $T_{0}$ space if and only if $M(x)=M(y)\Rightarrow
x=y.$ \ 

\item $(X,\tau )$ is (path and chain) connected if and only if for any pair
of points $a,b\in X$, there exists a finite set of points $%
\{a=x_{0},x_{1},...,x_{n}=b\}$ such that $N(x_{i-1})\cap N(x_{i})\neq
\varnothing $, $i\leq n.$
\end{enumerate}
\end{thm}

In some sense, a $T_{0}$ Alexandrov space is the most natural topological
structure induced by a partially ordered set. \ Indeed, if $(X,\leq )$ is
any partially ordered set, then the sets of the form $M(a)=\{x\in X:x\leq
a\} $ constitute a basis for a $T_{0}$ Alexandrov space $(X,\tau ).$ \
Conversely, if $(X,\tau )$ is a $T_{0}$ Alexandrov space, we can define a
relation $\leq $ on $X$ given by $a\leq b$ if and only if $a\in M(b).$ \
More precisely, we have the following result found in \cite{Arenas}

\begin{thm}
There is an isomorphism between the category of $T_{0}$ Alexandrov spaces
with continuous maps and the category of partially ordered sets with order
preserving set maps. \ 
\end{thm}

Now, let $(X,\tau )$ be any $T_{0}$ Alexandrov space with minimal
neighborhood base $\mathcal{M}=\{M(a):a\in X\}$. \ One can construct a
directed acyclic graph $\mathcal{G}(\mathcal{V},\mathcal{E})$ determined by
the space $(X,\tau ).$ \ The set $\mathcal{V}$ of vertices is taken to be
the underlying set $X$. \ Define an edge $a\rightarrow b$ if and only if $%
M(a)\subsetneq M(b)$ and there is no minimal base element $M(c)$ properly
between $M(a)$ and $M(b).$ \ The resulting construction is a simple (no
parallel edges) directed acyclic graph (SDAG). \ We now make the graphical
representation of the previous constructions precise.

\begin{defn}
Let $(X,\leq )$ be any partially ordered set and define intervals $(-\infty
,b]=\{x\in X:x\leq b\}$ and $[a,b]=\{x\in X:a\leq x\leq b\}$

\begin{enumerate}
\item We write $(X,\tau (\leq ))$ to denote the the Alexandrov topology
generated by the minimal neighborhood base $\mathcal{M}=\{M(a):a\in X\}$
where $M(a)=(-\infty ,a].$

\item We write $\mathcal{G}(X,\mathcal{E(\leq )})$ to denote the directed
acyclic graph whose vertices are the elements of $X$ and edges $a\rightarrow
b$ if and only if $a<b$ and $[a,b]=\{a,b\}.$
\end{enumerate}
\end{defn}

\begin{defn}
Let $\mathcal{V}$ be a nonempty set.

\begin{enumerate}
\item A \textit{finite directed path} is a sequence of edges $%
\{e_{1},e_{2},...,e_{n}\}\subset \mathcal{E}$ where $e_{i}=(v_{i-1},v_{i})$
for each $i\in \{1,2,...,n\}.$ \ A finite directed path in $\mathcal{G(V},%
\mathcal{E)}$ is also denoted by 
\begin{equation*}
v_{0}\rightarrow v_{1}\rightarrow v_{2}\rightarrow ...\rightarrow v_{n}.
\end{equation*}%
A directed graph is said to be \textit{acyclic} if there does not exist a
path $\{e_{1},e_{2},...,e_{n}\}\subset \mathcal{E}$ such that $v_{0}=v_{n}.$

\item A \textit{finite weak path} is a sequence of ordered pairs $%
\{e_{1},e_{2},...,e_{n}\}\subset \mathcal{V}\times \mathcal{V}$ where $%
e_{i}=(v_{i-1},v_{i})$ and either $(v_{i-1},v_{i})$ or $(v_{i},v_{i-1})\in 
\mathcal{E}.$ \ A finite directed path in $\mathcal{G(V},\mathcal{E)}$ is
also denoted by 
\begin{equation*}
v_{0}\leftrightarrow v_{1}\leftrightarrow ...\leftrightarrow v_{n}.
\end{equation*}

\item A directed graph $\mathcal{G(V},\mathcal{E)}$ is said to be \textit{%
weakly connected} if for every pair of vertices $v,w\in \mathcal{V}$ there
exists a finite weak path $\{e_{1},e_{2},...,e_{n}\}$ such that $v_{0}=v$
and $v_{n}=w.$
\end{enumerate}
\end{defn}

\section{The Graph of Divisibility}

In this section, we introduce the graph of divisibility of an integral
domain. \ We will see that this graph gives a picture of the group of
divisibility and can be used to detect certain factorization properties of a
domain. \ Although this graph does not detect all divisibility relations, it
does detect enough of the divisibility relation to clearly differentiate
atomicity and ACCP (for example). \ For the remainder of this article, we
denote the set of irreducible elements (atoms) of $D$ by $\text{Irr}(D)$
and the set of atomic elements (expressible as a finite product of atoms) is
denoted by $\mathcal{F}(D).$

\begin{defn}
Let $D$ be any integral domain with field of fractions $K$ and let $%
K^{\times }$ denote its multiplicative group$.$

\begin{enumerate}
\item We write $G(D)$ to denote the group of divisibility $K^{\times }/U(D)$
written additively. \ We write $G(D)^{+}$ to denote the positive elements of 
$G(D).$

\item We write $\mathcal{P}(D)$ to denote the group of nonzero principal
fractional ideals of $D$ partially ordered by inclusion. \ We write $%
\mathcal{P}(D)^{+}$ to denote the nonzero nonunit principal integral ideals
of $D.$
\end{enumerate}
\end{defn}

Recall that the ordering in $(G(D),\leq )$ is given by $\overline{a}\leq 
\overline{b}$ if and only if $\frac{b}{a}\in D.$ \ It readily follows that $%
0\leq \overline{a}$ if and only if $a\in D.$ \ It is easy to check that
there exists a reverse order group isomorphism $G(D)\rightarrow \mathcal{P}%
(D)$ given by $\overline{a}\mapsto aD$. \ With Definition 2.1 in hand, we
define a partial ordering on the set $\mathcal{P}(D)$ and consider the
structure of the associated topological space and directed acyclic graph. \
The following lemma is the basis for the remainder of our investigations.

\begin{lem}
Define a relation $\prec $\ on $\mathcal{P}$\ given by $a\prec b$\ if and
only if $\frac{a}{b}\in \mathcal{F}(D).$\ \ 

\begin{enumerate}
\item $(\mathcal{P}(D),\preceq )$\ is a partially ordered set.

\item $(\mathcal{P}(D),\tau (\preceq ))$\ is a $T_{0}$\ Alexandrov space
with neighborhood base given by the collection $\mathcal{M}(a)=\{x\in 
\mathcal{P}:\frac{x}{a}\in \mathcal{F}(D)\}.$

\item $\mathcal{G}(\mathcal{P}(D),\mathcal{E}(\preceq ))$\ is a directed
acyclic graph with directed edges $a\rightarrow b$\ if and only if $\frac{a}{%
b}\in \text{Irr}(D).$
\end{enumerate}
\end{lem}

\begin{proof}
(1) \ It is never the case that $a\prec a$ since $\frac{a}{a}$ is a unit,
and hence is not a product of atoms. \ Similarly, it is impossible that both 
$a\prec b$ and $a\succeq b$ can occur. \ Finally, if $a\prec b$ and $b\prec
c $, then $\frac{a}{b}\in \mathcal{F}(D)$ and $\frac{b}{c}\in \mathcal{F}%
(D). $ \ Since the set $\mathcal{F}(D)$ is multiplicatively closed, we have
that $\frac{a}{b}\cdot \frac{b}{c}=\frac{a}{c}\in \mathcal{F}(D)$ so that $%
a\prec c.$ \ 

\noindent (2) \ Follows immediately from (1) and the definition of $\prec .$
\ 

\noindent (3) \ If $a\rightarrow b$, then $a\prec b$ and $[a,b]=\{a,b\}$. \
It follows that $\frac{a}{b}=\pi _{1}\cdot \cdot \cdot \pi _{n}$ where each $%
\pi _{i}\in \text{Irr}(D)$. \ In other words, $[a,b]=\{a,\pi _{1}\cdot
\cdot \cdot \pi _{n-1}a,...,\pi _{1}a,b\}$ and the condition $[a,b]=\{a,b\}$
forces $b=\pi _{1}a.$ \ Therefore, $\frac{a}{b}\in \text{Irr}(D)$ as
needed. \ Conversely, if $\frac{a}{b}=\pi \in \text{Irr}(D)$, then it is
certainly true that $a\prec b$ and it suffices to check that $[a,b]=\{a,b\}.$
\ But if $a\prec x\prec b$, then $\frac{a}{x}=\pi _{1}\cdot \cdot \cdot \pi
_{n}$ and $\frac{x}{b}=\varsigma _{1}\cdot \cdot \cdot \varsigma _{m}$ where
each $\pi _{i},\varsigma _{i}\in \mathcal{A}(D).$ \ But then $\pi =\pi
_{1}\cdot \cdot \cdot \pi _{n}\varsigma _{1}\cdot \cdot \cdot \varsigma _{m}$
forcing (without loss of generality) $\pi =\pi _{1}$ with the remaining
factors units. \ It follows that $x=b$ as needed.
\end{proof}

With Lemma 3.2 in hand, we make the definition central to our study.

\begin{defn}
We call $\mathcal{G}(\mathcal{P}(D),\mathcal{E}(\preceq ))$ the \textit{%
graph of divisibility} of $D.$ \ We might also refer to the subgraph $%
\mathcal{G}(\mathcal{P}(D)^{+},\mathcal{E}(\preceq ))$ the graph of
divisibility.
\end{defn}

We illustrate this definition with a few easy examples.

\begin{ex}
\begin{enumerate}
\item Let $D$ be a one-dimensional Noetherian valuation domain. \ It is
well-known that $D$ is a PID with a unique nonzero prime ideal. \ So the the
elements of $\mathcal{P}(D)^{+}$ can be enumerated by the positive integers.
\ We write $\mathcal{P}(D)^{+}=\{\pi ,\pi ^{2},\pi ^{3},...\}$ where $\pi $
is a chosen generator of the unique maximal ideal. \ The graph of
divisibility $\mathcal{G}(\mathcal{P}(D),\mathcal{E}(\preceq ))$ is the
(branchless) tree that looks like 
\begin{equation*}
...\rightarrow \pi ^{2}\rightarrow \pi \rightarrow 1\rightarrow \frac{1}{\pi 
}\rightarrow \frac{1}{\pi ^{2}}\rightarrow ....
\end{equation*}%
Similarly, the subgraph $\mathcal{G}(\mathcal{P}(D)^{+},\mathcal{E}(\preceq
))$ looks like 
\begin{equation*}
...\rightarrow \pi ^{3}\rightarrow \pi ^{2}\rightarrow \pi
\end{equation*}

\item \ Let $D$ be a one-dimensional nondiscrete valuation domain. \ For the
sake of concreteness, we will say that the corresponding value group is $%
%TCIMACRO{\U{211a} }%
%BeginExpansion
\mathbb{Q}
%EndExpansion
$. \ In this example, there are no irreducible elements and hence no two
elements of $\mathcal{P}(D)$ are adjacent. \ It follows that $\mathcal{G}(%
\mathcal{P}(D),\mathcal{E}(\preceq ))$ is just the collection of vertices
corresponding to $\mathcal{P}(D)$ with no edges whatsoever. \ In fact, the
graph of divisibility of any antimatter domain (no irreducible elements)
consists of vertices only. \ Topologically speaking, $(\mathcal{P}(D),\tau
(\preceq ))$ is totally disconnected. \ That is, given any $a\in \mathcal{P}%
(D)$ we have that $M(a)=\{a\}.$ \ The same is certainly true for the
subspace $(\mathcal{P}(D)^{+},\tau (\preceq ))$ and the subgraph $\mathcal{G}%
(\mathcal{P}(D)^{+},\mathcal{E}(\preceq )).$
\end{enumerate}
\end{ex}

Recall that a \textit{sink} in a directed graph is a vertex with arrows in
but no arrows out. \ We have the following lemma.

\begin{lem}
Let $D$ be an integral domain and let $\mathcal{G}(P(D)^{+},\mathcal{E}%
(\preceq ))$ be the associated graph of divisibility. \ Then,an element $\pi
\in D^{\bullet }$ is irreducible in $D$ if and only if the node $\pi $ is a
sink in $\mathcal{G}(P(D)^{+},\mathcal{E}(\preceq )).$
\end{lem}

It is well known that $D$ is atomic if and only if every element of $%
G(D)^{+} $ can be written as a sum of minimal positive elements. \
Similarly, $D$ is satisfies ACCP if and only if every descending sequence of
elements in $G(D)^{+}$ stabilizes. \ As with the group of divisibility, the
graph of divisibility can be used to characterize the well-studied
factorization domains. \ We close this section with the following result.

\begin{thm}
Let $D$\ be an integral domain and let $\mathcal{G}(\mathcal{P}(D)^{+},%
\mathcal{E}(\preceq ))$\ be the associated graph of divisibility.

\begin{enumerate}
\item $D$\ is atomic if and only if for every non unit element $a\in 
\mathcal{P}(D)^{+}$, there exists a (finite) path originating from $a$\ that
terminates at an atom.

\item $D$\ satisfies ACCP if and only if for every $a\in \mathcal{P}(D)^{+}$%
, every path originating from $a$\ terminates at an atom.

\item $D$\ is a BFD if and only if for every $a\in \mathcal{P}(D)^{+}$,
every path originating from $a$\ terminates at an atom and there is an upper
bound on the lengths of all such paths.

\item $D$\ is an FFD if and only if for every $a\in \mathcal{P}(D)^{+}$,
every path originating from $a$\ terminates at an atom and there are
finitely many such paths.

\item $D$\ is an HFD if and only if for every $a\in \mathcal{P}(D)^{+}$,
every path originating from $a$\ terminates at an atom and all such paths
are of the same length.
\end{enumerate}
\end{thm}

\section{Some Connectedness Properties}

In this section, we consider the connectedness of the graph of divisibility.
\ To do this, we will examine the connectedness of the associated Alexandrov
topology. \ We conclude the section with a few examples.

\begin{thm}
Let $D$\ be an integral domain. \ The following statements for $a,b\in
K^{\times }$\ are equivalent.

\begin{enumerate}
\item There exist atoms $\pi _{i},\xi _{i}\in \text{Irr}(D)$ such that $%
\frac{a}{b}=\frac{\pi _{1}\cdot \cdot \cdot \pi _{n}}{\xi _{1}\cdot \cdot
\cdot \xi _{m}}.$

\item The points $a,b$\ belong to the same connected component in the
Alexandrov topology $(\mathcal{P}(D),\tau (\preceq )).$

\item There is a finite weak path connecting $a$\ to $b$\ in the graph of
divisibility $\mathcal{G}(\mathcal{P}(D),\mathcal{E}(\preceq )).$
\end{enumerate}
\end{thm}

\begin{proof}
(1)$\Rightarrow $(2) \ Suppose there exist atoms $\pi _{i},\xi _{i}\in 
\text{Irr}(D)$ such that $\frac{a}{b}=\frac{\pi _{1}\cdot \cdot \cdot \pi
_{n}}{\xi _{1}\cdot \cdot \cdot \xi _{m}}.$ \ To show that $a,b$ belong to
the same connected component, it suffices to show that $M(a)\cap M(b)$ is
nonempty. \ To this end, note that 
\begin{equation*}
a\xi _{1}\cdot \cdot \cdot \xi _{m}=c=b\pi _{1}\cdot \cdot \cdot \pi _{n}
\end{equation*}%
implies that $c\prec a$ because $\frac{c}{a}=\xi _{1}\cdot \cdot \cdot \xi
_{m}$. \ Similarly, we have that $c\prec b$, from which it follows that $%
c\in M(a)\cap M(b)$ as needed.

(2)$\Rightarrow $(3) \ If $a,b$ belong to the same connected component, then
there exists a finite set of points $\{a=x_{0},x_{1},...,x_{n}=b\}$ such
that $M(x_{i-1})\cap M(x_{i})\neq \varnothing $ for all $i\in \{1,2...,n\}.$
\ Hence, we can choose a $c_{i}\in M(x_{i-1})\cap M(x_{i})$ so that $\frac{%
c_{i}}{x_{i-1}},\frac{c_{i}}{x_{i}}\in \mathcal{F}(D)$, say $\frac{c_{i}}{%
x_{i-1}}=\xi _{1}\cdot \cdot \cdot \xi _{m}$ and $,\frac{c_{i}}{x_{i}}=\pi
_{1}\cdot \cdot \cdot \pi _{n}$ where $\pi _{i},\xi _{i}\in \text{Irr}%
(D). $ \ It follows that there are directed paths 
\begin{equation*}
c_{i}\rightarrow (x_{i-1}\xi _{1}\cdot \cdot \cdot \xi _{m-1})\rightarrow
...\rightarrow (x_{i-1}\xi _{1})\rightarrow x_{i-1}
\end{equation*}%
and 
\begin{equation*}
c_{i}\rightarrow (x_{i}\pi _{1}\cdot \cdot \cdot \pi _{n})\rightarrow
...\rightarrow (x_{i}\pi _{1})\rightarrow x_{i}.
\end{equation*}%
Hence, there is a weak path%
\begin{equation*}
x_{i-1}\leftarrow ...\leftarrow c_{i}\rightarrow ...\rightarrow x_{i}
\end{equation*}%
for all $i\in \{1,2...,n\}$, and so there is a weak path connecting $a$ to $%
b $.

(3)$\Rightarrow $(1) \ Suppose that $a,b$ are distinct points in the $T_{0}$
Alexandrov space $(\mathcal{P}(D),\mathcal{\tau }(\preceq ))$. Then there
exists a finite weak path connecting $a$ to $b$ say 
\begin{equation*}
a=x_{0}\leftrightarrow x_{1}\leftrightarrow ...\leftrightarrow x_{n}=b.
\end{equation*}%
Using induction on $n$, we suppose that the result is true for all $k<n.$ \
It follows from the existence of the weak path $a=x_{0}\leftrightarrow
x_{1}\leftrightarrow ...\leftrightarrow x_{n-1}$ that $\frac{a}{x_{n-1}}=%
\frac{\pi _{1}\cdot \cdot \cdot \pi _{n}}{\xi _{1}\cdot \cdot \cdot \xi _{m}}
$ where $\pi _{i},\xi _{i}\in \text{Irr}(D).$ \ Now observe that either $%
x_{n-1}\rightarrow b$ or $b\rightarrow x_{n-1}$. \ If $x_{n-1}\rightarrow b$
then by definition, we have $\frac{x_{n-1}}{b}=\pi \in \text{Irr}(D).$ \
It follows that $\frac{a}{b}=\frac{a}{x_{n-1}}\cdot \frac{x_{n-1}}{b}=\frac{%
\pi \pi _{1}\cdot \cdot \cdot \pi _{n}}{\xi _{1}\cdot \cdot \cdot \xi _{m}}$%
, and a similar argument handles the case.
\end{proof}

If $\mathfrak{F}$ is the subgroup of $K^{\times }$ generated by $\text{Irr%
}(D)$, then we can relate the number of connected components of $\mathcal{G}(%
\mathcal{P},\mathcal{E}(\preceq ))$ with the order of the quotient group $%
K^{\times }/\mathfrak{F}$ (a homomorphic image of the group of divisibility $%
K^{\times }/U(D)$). \ We immediately get the following result. \ 

\begin{cor}
There are 1-1 correspondences between the elements of $K^{\times }/\mathfrak{%
F}$, the connected components of $\mathcal{G}(\mathcal{P}(D),\mathcal{E}%
(\preceq ))$, and the connected components of $(\mathcal{P}(D),\tau (\preceq
)).$ \ 
\end{cor}

\begin{ex}
Consider the classical construction $D=%
%TCIMACRO{\U{2124} }%
%BeginExpansion
\mathbb{Z}
%EndExpansion
+x%
%TCIMACRO{\U{211a} }%
%BeginExpansion
\mathbb{Q}
%EndExpansion
\lbrack x].$ \ It is well-known that the irreducible elements of $D$ are the
primes $p\in 
%TCIMACRO{\U{2124} }%
%BeginExpansion
\mathbb{Z}
%EndExpansion
$ and $%
%TCIMACRO{\U{211a} }%
%BeginExpansion
\mathbb{Q}
%EndExpansion
\lbrack x]$-irreducible polynomials of the form $\pm 1+xq(x)$ where $q(x)\in 
%TCIMACRO{\U{211a} }%
%BeginExpansion
\mathbb{Q}
%EndExpansion
\lbrack x]$ (see \cite{CS})$.$ \ For each $a\in D$ let us write $%
a(x)=(a_{0},a_{1},a_{2},...)$ where $a_{0}\in 
%TCIMACRO{\U{2124} }%
%BeginExpansion
\mathbb{Z}
%EndExpansion
$ and $a_{i}\in 
%TCIMACRO{\U{211a} }%
%BeginExpansion
\mathbb{Q}
%EndExpansion
$ for all $i\geq 1.$ \ As with power series representations, we define the
order of $a$ to be the natural number $\text{ord}(a)=\min \{i\in 
%TCIMACRO{\U{2115} }%
%BeginExpansion
\mathbb{N}
%EndExpansion
:a_{i}\neq 0\}.$ \ It follows from \cite{CS} that $a(x)\in \mathcal{F}(D)$
if and only if $\text{ord}(a)=0.$ \ We will now show that two polynomials 
$a,b\in D$ belong to the same connected component of $(\mathcal{P}%
(D)^{+},\tau (\preceq ))$ if and only if $\text{ord}(a)=\text{ord}(b).$
\ Indeed, write $a(x)=x^{e_{0}}\overline{a}(x)$ and $b(x)=x^{f_{0}}\overline{%
b}(x)$ where $\text{ord}(\overline{a})=0=\text{ord}(\overline{b})$
(allowing $e_{0}=0=f_{0}$). $\ $If $a(x)$ is connected to $b(x)$, then by
there exist atoms $\pi _{i},\xi _{i}\in \text{Irr}(D)$ such that%
\begin{equation*}
\frac{a(x)}{b(x)}=\frac{\pi _{1}(x)\cdot \cdot \cdot \pi _{n}(x)}{\xi
_{1}(x)\cdot \cdot \cdot \xi _{m}(x)}.
\end{equation*}%
We now have the equation 
\begin{equation*}
x^{e_{0}}\overline{a}(x)\xi _{1}(x)\cdot \cdot \cdot \xi _{m}(x)=x^{f_{0}}%
\overline{b}(x)\pi _{1}(x)\cdot \cdot \cdot \pi _{n}(x)
\end{equation*}%
and one easily checks that 
\begin{equation*}
e_{0}=\text{ord}(x^{e_{0}})=\text{ord}(x^{e_{0}}\overline{a}\xi
_{1}\cdot \cdot \cdot \xi _{m})=\text{ord}(x^{f_{0}}\overline{b}\pi
_{1}\cdot \cdot \cdot \pi _{n})=\text{ord}(x^{f_{0}})=f_{0}
\end{equation*}%
For the converse, suppose that $e_{0}=f_{0}.$ \ Again, using the fact that $%
\text{ord}(\overline{a})=0=\text{ord}(\overline{b})$ is equivalent to $%
\overline{a},\overline{b}\in \mathcal{F}(D)$, we have the existence of atoms 
$\pi _{i},\xi _{i}\in \text{Irr}(D)$ such that 
\begin{equation*}
\frac{\overline{a}(x)}{\overline{b}(x)}=\frac{\pi _{1}(x)\cdot \cdot \cdot
\pi _{n}(x)}{\xi _{1}(x)\cdot \cdot \cdot \xi _{m}(x)}.
\end{equation*}%
On the other hand, $e_{0}=f_{0}$ implies 
\begin{equation*}
\frac{\overline{a}(x)}{\overline{b}(x)}=\frac{x^{e_{0}}\overline{a}(x)}{%
x^{f_{0}}\overline{b}(x)}=\frac{a(x)}{b(x)}.
\end{equation*}%
It follows that the distinct connected components of $(\mathcal{P}%
(D)^{+},\tau (\preceq ))$ are given by the set $\{\text{Irr}%
(D)=[2],[x],[x^{2}],...\}$. \ In other words, there is no weak path $%
x^{m}\leftrightarrow ...\leftrightarrow x^{n}$ in $\mathcal{G}(\mathcal{P}%
^{+},\mathcal{E}(\preceq ))$ whenever $m\neq n.$
\end{ex}

\begin{ex}
Let $x,y$ be indeterminates over the field $\mathbb{F}_{2}.$

\begin{enumerate}
\item Now let $X=\{x^{\alpha }:\alpha \in 
%TCIMACRO{\U{211a} }%
%BeginExpansion
\mathbb{Q}
%EndExpansion
^{+}\}$ and $Z_{1}=\{\frac{y^{k}}{x^{\alpha }}:\alpha \in 
%TCIMACRO{\U{211a} }%
%BeginExpansion
\mathbb{Q}
%EndExpansion
^{+},k\in 
%TCIMACRO{\U{2124} }%
%BeginExpansion
\mathbb{Z}
%EndExpansion
^{+},$ and $k\geq 2\}.$ \ We determine the number of connected components in
the graph of divisibility of the domain $D_{1}=\mathbb{F}%
_{2}[X,y,Z_{1}]_{(X,y,Z_{1})}.$ \ To do this, we first observe that the
integral closure of $D_{1}$ is the rank $2$ valuation domain $V=\mathbb{F}%
_{2}[X,Z]_{(X,Z)}$ where $Z=\{\frac{y}{x^{\alpha }}:\alpha \in 
%TCIMACRO{\U{211a} }%
%BeginExpansion
\mathbb{Q}
%EndExpansion
^{+}\}.$ \ The value group of $V$ is $%
%TCIMACRO{\U{2124} }%
%BeginExpansion
\mathbb{Z}
%EndExpansion
\oplus 
%TCIMACRO{\U{211a} }%
%BeginExpansion
\mathbb{Q}
%EndExpansion
$ ordered lexicographically and it is easy to check that every element of $%
V^{\bullet }$ is a unit multiple of $y^{k}$ or $x^{r}y^{k}$ where $(k,r)\in 
%TCIMACRO{\U{2124} }%
%BeginExpansion
\mathbb{Z}
%EndExpansion
^{+}\oplus 
%TCIMACRO{\U{211a} }%
%BeginExpansion
\mathbb{Q}
%EndExpansion
.$ \ It is not hard to check that every element of $\text{Irr}(D_{1})$
has value $(1,0).$ \ It is now an easy matter to check that the connected
components of $\mathcal{G}(\mathcal{P}(D_{1})^{+},\mathcal{E}(\preceq ))$
are given in terms of their values by%
\begin{equation*}
\{[(0,\alpha )]\}_{\alpha \in 
%TCIMACRO{\U{211a} }%
%BeginExpansion
\mathbb{Q}
%EndExpansion
^{+}}\cup \{[(k,0)]\}_{k\in 
%TCIMACRO{\U{2124} }%
%BeginExpansion
\mathbb{Z}
%EndExpansion
^{+}}\cup \{[(k,\alpha )]\}_{\substack{ k\geq 2  \\ \alpha <0}}.
\end{equation*}%
For example, consider the elements $f=x^{\frac{1}{2}}$ and $g=\frac{y^{3}}{%
x^{\frac{1}{3}}}$. \ Then $v(f)=(0,\frac{1}{2})$ and $v(g)=(3,\frac{1}{3}).$
\ Then $v(\frac{g}{f})=v(g)-v(f)=(3,-\frac{1}{6})$ cannot be written in the
form $m(1,0)$ where $m\in 
%TCIMACRO{\U{2124} }%
%BeginExpansion
\mathbb{Z}
%EndExpansion
.$ \ In other words, $\frac{g}{f}$ cannot be expressed as the quotient of
atomic elements.

\item If $Z_{2}=\{\frac{y^{k}}{x^{j}}:j\in 
%TCIMACRO{\U{2124} }%
%BeginExpansion
\mathbb{Z}
%EndExpansion
^{+},k\in 
%TCIMACRO{\U{2124} }%
%BeginExpansion
\mathbb{Z}
%EndExpansion
^{+},$ and $k\geq 2\}$ and $D_{2}=\mathbb{F}_{2}[x,y,Z_{2}]_{(x,y,Z_{2})}$,
then $(\mathcal{P}(D_{2})^{+},\tau (\preceq ))$ is a connected Alexandrov
space. \ Equivalently, the graph of divisibility $\mathcal{G}(\mathcal{P}%
(D_{2})^{+},\mathcal{E}(\preceq ))$ is weakly connected. \ One need only
check that the integral closure of $D_{2}$ has the discrete value group $%
%TCIMACRO{\U{2124} }%
%BeginExpansion
\mathbb{Z}
%EndExpansion
\oplus 
%TCIMACRO{\U{2124} }%
%BeginExpansion
\mathbb{Z}
%EndExpansion
$ ordered lexicographically. \ Again, it is not hard to check that every
element of $\text{Irr}(D_{2})$ has value $(0,1)$ or $(1,0)$ and given any 
$f,g\in D_{2}$, we have that $v(\frac{g}{f})=m(1,0)+n(0,1)$ where $m,n\in 
%TCIMACRO{\U{2124} }%
%BeginExpansion
\mathbb{Z}
%EndExpansion
.$
\end{enumerate}
\end{ex}

\section{Some Generalizations of Atomicity}

In this section, we show that a connected graph of divisibility gives rise
to a generalization of atomicity. \ 

\begin{defn}
Let $D$ be any integral domain.

\begin{enumerate}
\item $D$ is called \textit{almost atomic} if for every $a\in D^{\bullet }$,
there exist atoms $\{\pi _{i}\}\subset \text{Irr}(D)$ such that $a\pi
_{1}\cdot \cdot \cdot \pi _{n}\in \mathcal{F}(D).$

\item $D$ is called \textit{quasi atomic} if for every $a\in D^{\bullet }$,
there exists an element $b\in D$ such that $ab\in \mathcal{F}(D).$
\end{enumerate}
\end{defn}

It is easy to see that almost atomic implies quasi atomic. \ Also, if $D$ is
quasi atomic, it is not hard to show that every nonzero prime ideal of $D$
contains an irreducible element. \ We have the following lemma.

\begin{lem}
Given an integral domain $D$, each condition below implies the next:

\begin{enumerate}
\item $D$ is atomic

\item $D$ is almost atomic

\item $D$ is quasi atomic

\item Every nonzero prime ideal of $D$ contains an irreducible element.
\end{enumerate}
\end{lem}

\begin{proof}
It suffices to show that (3) implies (4). \ Suppose that $D$ is quasi
atomic. \ If $a$ is a nonzero element of a prime ideal $P$, then there is $%
b\in D$ such that $ab=\pi _{1}\cdot \cdot \cdot \pi _{n}$ where each $\pi
_{i}\in \text{Irr}(D).$ \ But then $\pi _{1}\cdot \cdot \cdot \pi _{n}\in
P$ so that $\pi _{i}\in P$ for some $i\leq n.$ \ 
\end{proof}

These observations give an example of an integral domain that is not quasi
atomic.

\begin{ex}
As in Example 4.3, let $D=%
%TCIMACRO{\U{2124} }%
%BeginExpansion
\mathbb{Z}
%EndExpansion
+x%
%TCIMACRO{\U{211a} }%
%BeginExpansion
\mathbb{Q}
%EndExpansion
\lbrack x].$ \ Then $x%
%TCIMACRO{\U{211a} }%
%BeginExpansion
\mathbb{Q}
%EndExpansion
\lbrack x]$ is a prime ideal of $D$ that contains no irreducible element. \
To see this, recall from 4.3 that if $f\in \text{Irr}(D)$, then $\text{%
ord}(f)=0.$ \ But $f\in x%
%TCIMACRO{\U{211a} }%
%BeginExpansion
\mathbb{Q}
%EndExpansion
\lbrack x]$ if and only if $\text{ord}(f)\geq 1.$ \ It follows from Lemma
5.2 that $D$ is not quasi atomic.
\end{ex}

We now show the connection between almost atomicity and a connected graph of
divisibility.

\begin{thm}
The following statements are equivalent for a domain $D$.

\begin{enumerate}
\item $D$\ is almost atomic.

\item $(\mathcal{P},\tau (\preceq ))$\ is connected.

\item $\mathcal{G}(\mathcal{P},\mathcal{E}(\preceq ))$ is weakly connected.
\end{enumerate}
\end{thm}

\begin{proof}
(1)$\Rightarrow $(2) \ Choose any two points $a,b\in (\mathcal{P},\tau
(\preceq ))$. \ If $D$ is almost atomic, there exist atoms $\pi _{i},\xi
_{i}\in \text{Irr}(D)$ such that $a\pi _{1}\cdot \cdot \cdot \pi _{n}$
and $b\xi _{1}\cdot \cdot \cdot \xi _{m}\in \mathcal{F}(D).$ \ In other
words, there exist atoms $\sigma _{i},\varsigma _{i}\in \mathcal{A}(D)$ such
that $\frac{a}{b}=\frac{\pi _{1}\cdot \cdot \cdot \pi _{n}\sigma _{i}\cdot
\cdot \cdot \sigma _{j}}{\xi _{1}\cdot \cdot \cdot \xi _{m}\varsigma
_{i}\cdot \cdot \cdot \varsigma _{k}}.$ \ It follows from Theorem 4.1 that
any pair of points in $(\mathcal{P},\tau (\preceq ))$ belong to the same
connected component. \ 

(2)$\Rightarrow $(3) \ Follows immediately from Theorem 4.1.

(3)$\Rightarrow $(1) \ Since $\mathcal{G}(\mathcal{P},\mathcal{E}(\preceq ))$
is weakly connected, there is a weak path connecting any $\frac{a}{1}\in 
\mathcal{P}$ (where $a\in D$) to an element of the form $\frac{\pi }{1}$
where $\pi \in \text{Irr}(D).$ \ Theorem 4.1 implies that there exist
atoms $\pi _{i},\xi _{i}\in \text{Irr}(D)$ such that $\frac{a}{\pi }=%
\frac{\pi _{1}\cdot \cdot \cdot \pi _{n}}{\xi _{1}\cdot \cdot \cdot \xi _{m}}%
.$ \ In other words, there exist atoms $\xi _{i}\in \text{Irr}(D)$ such
that $a\xi _{1}\cdot \cdot \cdot \xi _{m}\in \mathcal{F}(D).$
\end{proof}

Using the Theorem 5.4 and the results from the previous section, we are led
to an example of an almost atomic domain that is not atomic.

\begin{ex}
\begin{enumerate}
\item As in Example 4.4(1), let $D_{1}=\mathbb{F}%
_{2}[X,y,Z_{1}]_{(X,y,Z_{1})}$ where $X=\{x^{\alpha }:\alpha \in 
%TCIMACRO{\U{211a} }%
%BeginExpansion
\mathbb{Q}
%EndExpansion
^{+}\}$ and $Z_{1}=\{\frac{y^{k}}{x^{\alpha }}:\alpha \in 
%TCIMACRO{\U{211a} }%
%BeginExpansion
\mathbb{Q}
%EndExpansion
^{+},k\in 
%TCIMACRO{\U{2124} }%
%BeginExpansion
\mathbb{Z}
%EndExpansion
^{+},$ and $k\geq 2\}.$ \ Since the connected components are in a 1-1
correspondence with $%
%TCIMACRO{\U{211a} }%
%BeginExpansion
\mathbb{Q}
%EndExpansion
$, it is certainly \textit{not} the case that $D_{1}$ is almost atomic
(Theorem 5.4). \ However, it is the case that $D_{1}$ is quasi atomic. \
Indeed, given any $f\in D_{1}^{\bullet }$, we can write $v(f)=(k,\alpha )$.
\ There is a $g\in D_{1}$ such that $v(g)=(2,-\alpha )$ and so 
\begin{equation*}
v(fg)=v(f)+v(g)=(k+2,0)=(1,0)+...+(1,0).
\end{equation*}%
Translating this information back to $D_{1},$ we get that $fg=y^{k+2}u$ for
some unit $u\in V.$ \ Note that if $y^{n+1}u\in D_{1}$ for some unit $u\in V$%
, then $y^{n+1}u=v_{1}x^{\alpha }+v_{2}\frac{y^{l}}{x^{\beta }}+v_{3}y$
where either $v_{i}\in U(D_{1})$ or $v_{i}=0.$ \ If $n>0$, then $%
v_{1}=0=v_{3}.$ \ It follows that $l\geq n+1$ so that $y^{n}u=v_{2}\frac{%
y^{l-1}}{x^{\beta }}.$ \ Therefore, $yu=v_{2}\frac{y^{l-n}}{x^{\beta }}\in
D_{1}$ as $l-n\geq 1.$ \ It follows from all of this that $fg\in \mathcal{F}%
(D)$ as needed.

\item As in Example 4.4(2), let $D_{2}=\mathbb{F}%
_{2}[x,y,Z_{2}]_{(x,y,Z_{2})}$ where $Z_{2}=\{\frac{y^{k}}{x^{j}}:j\in 
%TCIMACRO{\U{2124} }%
%BeginExpansion
\mathbb{Z}
%EndExpansion
^{+},k\in 
%TCIMACRO{\U{2124} }%
%BeginExpansion
\mathbb{Z}
%EndExpansion
^{+},$ and $k\geq 2\}.$ \ Since $\mathcal{G}(\mathcal{P}(D_{2})^{+},\mathcal{%
E}(\preceq ))$ is weakly connected, it must be the case that $D_{2}$ is
almost atomic. \ However, it is not atomic since, for example, $v(\frac{y^{2}%
}{x^{\frac{1}{2}}})=(2,\frac{1}{2})$ cannot be written as an $%
%TCIMACRO{\U{2115} }%
%BeginExpansion
\mathbb{N}
%EndExpansion
$-linear combination $m(1,0)+n(0,1).$
\end{enumerate}
\end{ex}

\subsection*{Acknowledgements}

The authors would like to thank the North Dakota State University Department
of Mathematics for their continued support.

\bibliography{biblio}{}
\bibliographystyle{plain}

\end{document}